\documentclass[11pt,reqno]{amsart}
\usepackage{amsthm}
\usepackage{amssymb}
\usepackage{latexsym}
\usepackage{multicol}
\usepackage{verbatim,enumerate}
\usepackage{amscd}
\usepackage{hyperref}
\usepackage{amsmath, amscd}
\usepackage{color}

\advance\textwidth by 1.2in \advance\oddsidemargin by -.6in \advance\evensidemargin by -.6in
\newtheorem*{cor}{Corollary}
\newtheorem*{lem}{Lemma}
\newtheorem*{prop}{Proposition}

\theoremstyle{definition}
\newtheorem*{defn}{Definition}
\theoremstyle{definition}
\newtheorem{thm}{Theorem}
\newtheorem*{conj}{Conjecture}
\newtheorem*{rem}{Remark}

\newenvironment{pf}{\proof}{\endproof}
\newcounter{cnt}
 \makeatletter
\def\mydggeometry{\makeatletter\dg@YGRID=1\dg@XGRID=20\unitlength=0.003pt\makeatother}
\makeatother \theoremstyle{remark}

\numberwithin{equation}{section}

 \DeclareMathOperator{\Ht}{ht} 
\let\bwdg\bigwedge
\def\bigwedge{{\textstyle\bwdg}}

\DeclareMathOperator{\Max}{Max}
\DeclareMathOperator{\ann}{Ann}
\DeclareMathOperator{\sym}{sym}

\begin{document}

\newcommand{\thmref}[1]{Theorem~\ref{#1}}
\newcommand{\secref}[1]{Section~\ref{#1}}
\newcommand{\lemref}[1]{Lemma~\ref{#1}}
\newcommand{\propref}[1]{Proposition~\ref{#1}}
\newcommand{\corref}[1]{Corollary~\ref{#1}}
\newcommand{\remref}[1]{Remark~\ref{#1}}
\newcommand{\defref}[1]{Definition~\ref{#1}}
\newcommand{\er}[1]{(\ref{#1})}
\newcommand{\id}{\operatorname{id}}
\newcommand{\ord}{\operatorname{\emph{ord}}}
\newcommand{\sgn}{\operatorname{sgn}}
\newcommand{\wt}{\operatorname{wt}}
\newcommand{\tensor}{\otimes}
\newcommand{\from}{\leftarrow}
\newcommand{\nc}{\newcommand}
\newcommand{\rnc}{\renewcommand}
\newcommand{\dist}{\operatorname{dist}}
\newcommand{\qbinom}[2]{\genfrac[]{0pt}0{#1}{#2}}
\nc{\cal}{\mathcal} \nc{\goth}{\mathfrak} \rnc{\bold}{\mathbf}
\renewcommand{\frak}{\mathfrak}
\newcommand{\supp}{\operatorname{supp}}
\renewcommand{\Bbb}{\mathbb}
\nc\bomega{{\mbox{\boldmath $\omega$}}} \nc\bpsi{{\mbox{\boldmath $\Psi$}}}
 \nc\balpha{{\mbox{\boldmath $\alpha$}}}
 \nc\bpi{{\mbox{\boldmath $\pi$}}}
\nc\bsigma{{\mbox{\boldmath $\sigma$}}} \nc\bcN{{\mbox{\boldmath $\cal{N}$}}} \nc\bcm{{\mbox{\boldmath $\cal{M}$}}} \nc\bLambda{{\mbox{\boldmath
$\Lambda$}}}

\newcommand{\lie}[1]{\mathfrak{#1}}
\makeatletter
\def\section{\def\@secnumfont{\mdseries}\@startsection{section}{1}%
  \z@{.7\linespacing\@plus\linespacing}{.5\linespacing}%
  {\normalfont\scshape\centering}}
\def\subsection{\def\@secnumfont{\bfseries}\@startsection{subsection}{2}%
  {\parindent}{.5\linespacing\@plus.7\linespacing}{-.5em}%
  {\normalfont\bfseries}}
\makeatother
\def\subl#1{\subsection{}\label{#1}}
 \nc{\Hom}{\operatorname{Hom}}
  \nc{\mode}{\operatorname{mod}}
\nc{\End}{\operatorname{End}} \nc{\wh}[1]{\widehat{#1}} \nc{\Ext}{\operatorname{Ext}} \nc{\ch}{\text{ch}} \nc{\ev}{\operatorname{ev}}
\nc{\Ob}{\operatorname{Ob}} \nc{\soc}{\operatorname{soc}} \nc{\rad}{\operatorname{rad}} \nc{\head}{\operatorname{head}}

 \nc{\Cal}{\cal} \nc{\Xp}[1]{X^+(#1)} \nc{\Xm}[1]{X^-(#1)}
\nc{\on}{\operatorname} \nc{\Z}{{\bold Z}} \nc{\J}{{\cal J}} \nc{\C}{{\bold C}} \nc{\Q}{{\bold Q}}
\renewcommand{\P}{{\cal P}}
\nc{\N}{{\Bbb N}} \nc\boa{\bold a} \nc\bob{\bold b} \nc\boc{\bold c} \nc\bod{\bold d} \nc\boe{\bold e} \nc\bof{\bold f} \nc\bog{\bold g}
\nc\boh{\bold h} \nc\boi{\bold i} \nc\boj{\bold j} \nc\bok{\bold k} \nc\bol{\bold l} \nc\bom{\bold m} \nc\bon{\bold n} \nc\boo{\bold o}
\nc\bop{\bold p} \nc\boq{\bold q} \nc\bor{\bold r} \nc\bos{\bold s} \nc\boT{\bold t} \nc\boF{\bold F} \nc\bou{\bold u} \nc\bov{\bold v}
\nc\bow{\bold w} \nc\boz{\bold z} \nc\boy{\bold y} \nc\ba{\bold A} \nc\bb{\bold B} \nc\bc{\bold C} \nc\bd{\bold D} \nc\be{\bold E} \nc\bg{\bold
G} \nc\bh{\bold H} \nc\bi{\bold I} \nc\bj{\bold J} \nc\bk{\bold K} \nc\bl{\bold L} \nc\bm{\bold M} \nc\bn{\bold N} \nc\bo{\bold O} \nc\bp{\bold
P} \nc\bq{\bold Q} \nc\br{\bold R} \nc\bs{\bold S} \nc\bt{\bold T} \nc\bu{\bold U} \nc\bv{\bold V} \nc\bw{\bold W} \nc\bz{\bold Z} \nc\bx{\bold
x} \nc\KR{\bold{KR}} \nc\rk{\bold{rk}} \nc\het{\text{ht }}

\nc\toa{\tilde a} \nc\tob{\tilde b} \nc\toc{\tilde c} \nc\tod{\tilde d} \nc\toe{\tilde e} \nc\tof{\tilde f} \nc\tog{\tilde g} \nc\toh{\tilde h}
\nc\toi{\tilde i} \nc\toj{\tilde j} \nc\tok{\tilde k} \nc\tol{\tilde l} \nc\tom{\tilde m} \nc\ton{\tilde n} \nc\too{\tilde o} \nc\toq{\tilde q}
\nc\tor{\tilde r} \nc\tos{\tilde s} \nc\toT{\tilde t} \nc\tou{\tilde u} \nc\tov{\tilde v} \nc\tow{\tilde w} \nc\toz{\tilde z} \nc\woi{w_{\omega_i}}

\author[M. Bennett, V. Chari, J. Greenstein, N. Manning]{Matthew Bennett, Vyjayanthi Chari, Jacob Greenstein, Nathan Manning}\address{ Department of Mathematics, University of California, Riverside, CA 92521}
\email{mbenn002@math.ucr.edu}
\email{vyjayanthi.chari@ucr.edu}\email{jacob.greenstein@ucr.edu}
\email{nmanning@math.ucr.edu}

\thanks{Partially supported by DMS-0901253 (V.~C.) and DMS-0654421 (J.~G.)}

\date{\today}

\title{On homomorphisms between global Weyl modules}

\begin{abstract} Let $\lie g$ be a
simple finite dimensional Lie algebra and let $A$ be a commutative associative algebra with unity.
Global Weyl modules for the generalized loop algebra $\lie  g\tensor A$ were defined in \cite{CPweyl,FL} for any dominant integral weight $\lambda$ of $\lie g$ by generators and
relations and further studied in \cite{CFK}. They are expected to play a role similar to that of
Verma modules in the study of categories of representations of $\lie g\tensor A$.
One of the fundamental properties of Verma modules is that the space of  morphisms between
two Verma modules is either zero or one--dimensional and also that any non--zero morphism is  injective. The aim of this paper is to
establish an analogue of this property for global Weyl modules. This is done under
certain restrictions on $\lie g$, $\lambda$ and $A$. A crucial tool is the
construction of fundamental global Weyl modules in terms of fundamental local
Weyl modules given in Section 3.

\end{abstract}
\maketitle
\section*{Introduction} In this paper, we continue the study of   representations of loop algebras and, more generally, of Lie algebras of the form $\lie g\otimes A$, where  $\lie g$ is a complex finite--dimensional simple Lie algebra and $A$ is a commutative associative algebra with unity over the complex numbers. To be precise, we are interested in the category $\mathcal I_A$ of   $\lie g\otimes A$--modules which are integrable as $\lie g$--modules. One motivation for studying this category is that it is closely related to the  categories
$\mathcal F$ and $\mathcal F_q$ of finite--dimensional representations, respectively, of affine Lie and quantum affine algebras, which has been a subject of considerable interest in recent years. The categories $\cal I_A$, $\mathcal F$ and $\cal F_q$ fail to be semi-simple, and  it was proved in \cite{CPweyl}  that irreducible representations of the quantum affine algebra specialize at $q=1$ to reducible indecomposable  representations of the loop algebra (obtained by taking $A=\bc[t,t^{-1}]$). This phenomenon is analogous to the one observed in modular representation theory, where an irreducible finite--dimensional  representation in characteristic zero becomes reducible by passing to characteristic $p$ and is called a Weyl module.

This analogy motivated the definition of  Weyl modules (global and local) for loop algebras in~\cite{CPweyl}. Their study was pursued for more
general rings~$A$ in~\cite{CFK} and~\cite{FL}.  Thus, given any dominant integral weight of the semi-simple Lie algebra $\lie g$
one can define an infinite--dimensional object
$W_A(\lambda)$ of~$\cal I_A$, called the
{\em global Weyl module}, via generators and relations. It was shown (see \cite{CFK} for the most general case) that if $A$ is finitely generated, then   $W_A(\lambda)$ is  a  right module for a certain commutative  finitely generated associative algebra $\ba_\lambda$, which is canonically associated with $A$ and $\lambda$. The local Weyl modules are obtained by tensoring the global Weyl module $W_A(\lambda)$ over $\ba_\lambda$ with simple $\ba_\lambda$-modules or, equivalently, can be given via generators and relations.

The local Weyl modules have been useful in understanding the blocks of the category $\cal F_A$  of finite--dimensional representations of $\lie g\otimes A$, and the quantum analogs are used to understand the blocks in $\cal F_q$. One motivation for this paper is to explore the use of the global Weyl modules to further understand the homological properties of the  categories $\cal F_A$ and $\cal I_A$. The global Weyl modules have nice universal
properties, and one expects them to play a role similar to that of the Verma modules $M(\lambda)$ in the study of the Bernstein-Gelfand-Gelfand category $\cal O$ for $\lie g$.  One of the most basic results about Verma modules is that $\Hom_{\lie g}(M(\lambda), M(\mu))$ is of dimension at most one and also that any non--zero map is injective. In this paper we prove the following analogue of this result.
\begin{thm}
Suppose that $\lie g=\lie{sl}_{r+1}$ or~$\lie g=\lie{sp}_{2r}$ and $A=\bc[t]$ or~$\bc[t^{\pm1}]$. Then
$$
\Hom_{\lie g\tensor A}(W_A(\mu),W_A(\lambda))\cong_{\ba_\lambda} \begin{cases}
0,& \mu\not=\lambda\\
\ba_\lambda,&\mu=\lambda
\end{cases}
$$
and every non-zero element of $\Hom_{\lie g\tensor A}(W_A(\lambda),W_A(\lambda))$ is injective.
\end{thm}
A more general version of this statement is provided in Theorem~\ref{thm2}.
Note that instead of having dimension one, the $\Hom$-spaces become free modules of rank one over~$\ba_\lambda$. This of course reflects the
fact that the ``top'' weight space of $M(\lambda)$ is one-dimensional, while for $W_A(\lambda)$ it is isomorphic to $\ba_\lambda$ as
a $\ba_\lambda$-module. It should be noted that, although the vanishing of $\Hom$ for $\mu\not=\lambda$ definitely fails for
$\lie g$ of other types, as shown in the remark in~\S\ref{red}, we still
expect that all non-zero elements of $\Hom_{\lie g\tensor A}(W_A(\mu),W_A(\lambda))$ are injective.

As is natural, one is guided by the representation theory of quantum
affine algebras and the phenomena occurring in $\cal F_q$. It is
interesting to note that Theorem~\ref{thm2}, for instance,  fails
exactly when the irreducible finite--dimensional module in the
quantum case specializes to a reducible module for the loop algebra.
We conclude by noting that global Weyl modules are
 also  defined for quantum affine algebras and similar questions can
 be posed in the quantum case. Some of these have have been studied
  in~\cite{BN} using crystal bases. However, it should be noted that
the results in the quantum case do not specialize to the classical
case. The results in the classical case are sometimes different from
those in the quantum situation. One of the reasons is that the local
fundamental Weyl module is always irreducible in the quantum case,
while its classical analogue is often reducible.

We now explain in more detail the organization and  results of this
paper. In Section~\ref{sect:loc glob wey} we fix the notation and
recall the basic facts on global Weyl modules that will be necessary
for the sequel. The main results are formulated and discussed in
Section~\ref{sect:main res}, together with their relation with the
quantum case. Section~\ref{sect:local Weyl} contains properties of
the algebras $\ba_\lambda$ and local Weyl modules which are needed
later. One of the principal difficulties in working with global Weyl
modules is their  abstract definition. However, when $A$ is the ring
of (generalized) Laurent polynomials, it admits a natural bialgebra
structure which can be used to reconstruct a {\em fundamental}
global Weyl module $W_A(\omega_i)$ from a local Weyl module. This
explicit realization of fundamental Weyl modules plays a crucial
role in proving our main results and occupies Section~\ref{sect:Loc
and glob Weyl} (Proposition~\ref{glfundreal}). It is not hard to
deduce from the definition of global Weyl modules that there exists
a canonical map from $W_A(\lambda)$ to a suitable tensor product of
global fundamental Weyl modules. Our next result (Theorem \ref{thm1}
and its Corollary~\ref{class}) studies the more general problem of
describing the space of $\lie g\tensor A$-module homomorphisms from
$W_A(\lambda)$ to an arbitrary tensor product of fundamental global
Weyl modules. The structure of the fundamental local Weyl modules
was given in \cite{CFK} for all classical simple Lie algebras and
for certain nodes of exceptional Lie algebras. This information and
the realziation of fundamental global Weyl modules provided by
Proposition~\ref{glfundreal} allows one to establish Theorem
\ref{thm1}.

In Theorem~\ref{thm2}, we use Theorem \ref{thm1} to compute
 $\Hom_{\lie g\otimes A}(W_A(\lambda), W_A(\mu))$ when $A$ is the ring of polynomials or Laurent polynomials in one  variable and under suitable conditions on $\mu$. We are able to do this because one has a lot of information (see
\cite{BN,CL,CPweyl,FoL,Ka,N})
on global Weyl modules in this case. We are also able to use a result of~\cite{FL} and Theorem~\ref{thm1} of this paper to compute $\Hom_{\lie g\otimes A}(W_A(\lambda), W_A(\mu))$ when $\lie g$ is of type $\lie{sl}_{r+1}$,  $\mu$ is a multiple of the first fundamental weight and $A$ is the (generalized) ring of Laurent polynomials in several  variables. This is done in Section~\ref{sect:pf thm2} of the paper.


\section{Global  Weyl modules}\label{sect:loc glob wey}
In this section we establish the notation to be
used in the rest of the paper and then recall the definition and
some elementary properties of the  global  Weyl modules.
\medskip

\subsection{}\label{subs:prelim.10}
Let  $\bc$ be the field of complex numbers and let $\bz$ (respectively $\bz_+$) be the set of integers (respectively non--negative integers). Given two complex vector spaces $V$, $W$ let $V\otimes W$ (respectively, $\Hom(V,W)$) denote their  tensor product over $\bc$ and (respectively the space of  $\bc$--linear maps from $V$ to $W$).

Given a commutative and  associative algebra $A$ over $\bc$, let  $\Max A$ be the  maximal spectrum of $A$ and $\mod A$ the category of left $A$--modules. Given
a right $A$-module $M$ and an element
$m\in M$, the annihilating (right) ideal of $m$ is    $$\ann_A m=\{a\in A\,:\, m.a = 0\}.$$
\subsection{}\label{subs:prelim.20}
For a complex Lie algebra $\lie a$  let $\bu(\lie a)$ be the associated
universal enveloping algebra. It is that the assignment $x\mapsto
x\otimes 1+1\otimes x$ for $x\in\lie a$ extends to a homomorphism
of algebras $\Delta:\bu(\lie a)\to \bu(\lie a)\otimes \bu(\lie a)$ and
defines a bialgebra structure on~$\bu(\lie a)$.
In particular, if
$V,W$ are two $\lie a$--modules then $V\tensor W$  and $\Hom_\bc(V,W)$ are naturally  $\bu(\lie a)$-modules
and $W\tensor V\cong V\tensor W$ as $\bu(\lie a)$-modules. One can also define the trivial $\lie a$--module structure on $\bc$ and we have
$$V^{\lie a}=\{v\in V: \lie a v=0\}\cong\Hom_{\lie a}(\bc, V).$$

Suppose that $A$ is an associative commutative algebra over $\bc$
with unity.  Then $\lie a\otimes A$ is canonically a Lie algebra,
with the Lie bracket given by
$$[x\otimes a, y\otimes b]=[x,y]\otimes ab,\qquad x,y\in\lie a,\qquad a,b\in A.
$$
We shall identify $\lie a$  with the
Lie subalgebra $\lie a\otimes 1$ of $\lie a\otimes A$. Note that for any algebra homomorphism  $\varphi:A\to A'$ 
the canonical map $1\tensor \varphi:\lie g\tensor A\to \lie g\tensor A'$ is a homomorphism of Lie algebras and hence induces an algebra homomorphism  $\bu(\lie g\tensor A)\to \bu(\lie g\tensor A')$.

\subsection{}
\label{C:section 2.1} Let $\lie g$ be a finite--dimensional
 complex simple Lie algebra with  a fixed  Cartan subalgebra $\lie h$.   Let $\Phi$  be the
  corresponding root system and  fix a set $\{\alpha_i: i\in I\}\subset\lie h^* $ (where $I=\{1,\dots,\dim\lie h\}$) of simple roots for $\Phi$.
   The root lattice $Q$ is the $\bz$--span of the simple roots while $Q^+$  is the $\bz_+$--span of the simple roots, and $\Phi^+=\Phi\cap  Q^+$ denotes the set of positive roots in $\Phi$. Let  $\Ht:Q^+\to\bz_+$ be the  homomorphism of free semi-groups defined by setting $\Ht(\alpha_i)=1$, $i\in I$.

   The restriction of  the Killing  form $\kappa:\lie g\times\lie g\to\bc$   to $\lie h\times\lie h$ induces a non--degenerate bilinear form $(\cdot ,\cdot )$ on $\lie h^*$,  and we let $\{\omega_i\,:\,i\in I\}\subset\lie h^*$ be the fundamental weights defined
by $2(\omega_j,\alpha_i)=\delta_{i,j}(\alpha_i,\alpha_i)$, $i,j\in I$.
 Let $P$ (respectively $P^+$) be the $\bz$ (respectively $\bz_+$) span of the $\{\omega_i\,:\,i\in I\}$ and note that $Q\subseteq P$. Given $\lambda,\mu\in P$ we say that $\mu\le \lambda$ if and only if $\lambda-\mu\in Q^+$.
 Clearly $\le $ is a partial order on $P$. The set $\Phi^+$ has a unique maximal element with respect to this order which is denoted by $\theta$ and is called the highest root of $\Phi^+$.  From now on, we normalize the bilinear form on $\lie h^*$ so that $(\theta, \theta)=2$.

 \subsection{} Given $\alpha\in \Phi$, let $\lie g_\alpha$ be the corresponding root space and define subalgebras
 $\lie n^\pm$ of $\lie g$ by $$\lie n^\pm=\bigoplus_{\alpha\in \Phi^+}\lie g_{\pm\alpha}.$$
 We have isomorphisms of vector spaces
  \begin{equation}\label{C:gtriang}\lie g\cong\lie n^-\oplus\lie h\oplus\lie n^+,\ \ \  \bu(\lie
g)\cong\bu(\lie n^-)\otimes\bu(\lie h)\otimes \bu(\lie
n^+). \end{equation}
   For $\alpha\in \Phi^+$, fix elements $x^\pm_\alpha\in\lie g_{\pm\alpha}$ and $h_\alpha\in\lie h$ spanning a Lie subalgebra of $\lie g$ isomorphic to $\lie{sl}_2$, i.e., we have $$[h_\alpha,x^\pm_\alpha]=\pm 2 x^\pm_\alpha,\qquad [x^+_\alpha, x^-_\alpha]= h_\alpha, $$ and more generally, assume that the set $\{x^\pm_\alpha: \alpha\in\Phi^+\}\cup\{h_i:=h_{\alpha_i}\,:\, i\in I\}$ is a Chevalley basis for~$\lie g$.

\subsection{} Given an $\lie h$--module $V$, we say that $V$ is a weight
module\index{weight module} if $$V=\bigoplus_{\mu\in\lie h^*}
V_\mu,\ \ V_\mu=\{v\in V: hv=\mu( h)v,\, h\in\lie h\} ,$$ and we set $\wt V=\{\mu\in\lie h^*: V_\mu\ne 0\}$.
If $\dim V_\mu<\infty$ for alll $\mu\in\lie h^*$, let $\ch\ V$ be the character of $V$, namely the element of the group ring $\bz[\lie h^*]$ given by, $$ \ch\  V=\sum_{\mu\in\lie h^*}\dim V_\mu e(\mu),$$ where  $e(\mu)\in\bz[\lie h^*]$ is the element  corresponding to~$\mu\in\lie h^*$. Observe that for two such modules $V_1$ and $V_2$, we have $$\ch(V_1\oplus V_2)=\ch \ V_1+\ch\ V_2,\qquad \ch  (V_1\otimes V_2)=\ch\ V_1 \ch\ V_2. $$

  \subsection{}  For $\lambda\in P^+$,
let $V(\lambda)$ be the left $\lie g$--module  generated by an
element $v_\lambda$  with defining
relations:$$hv_\lambda=\lambda(h)v_\lambda,\qquad
x_{\alpha_i}^+v_\lambda=0,\qquad
(x_{\alpha_i}^-)^{\lambda(h_{\alpha_i})+1} v_\lambda=0,$$\ where
$h\in\lie h$ and $i\in I$. It is well--known that $V(\lambda)$ is an
irreducible weight module, and  $$\dim V(\lambda)<\infty,\qquad \dim
V(\lambda)_\lambda=1,\qquad \wt V(\lambda)\subset\lambda-Q^+.$$ If
$V$ is  a finite-dimensional irreducible $\lie g$--module then there
exists a unique $\lambda\in P^+$ such that $V$ isomorphic to to
$V(\lambda)$. We shall say that  $V$ is a locally
finite--dimensional $\lie g$--module if $$\dim\bu(\lie g)v<\infty,\
\ v\in V.$$ It is well--known that a locally  finite--dimensional
$\lie g$--module is isomorphic to a direct sum of irreducible
finite--dimensional modules, moreover $$V_\mu\cap V^{\lie
n^+}\cong\Hom_{\lie g}(V(\mu), V),\ \ \mu\in P^+.$$

\subsection{}\label{defn-gl-weyl} Assume from now on that $A$ is an associative commutative algebra over $\bc$ with
unity. We recall the definition of the global Weyl modules. These
were first introduced and studied in the case when $A=\bc[t,t^{-1}]$
in \cite{CPweyl} and then later in \cite{FL} in the general case.
 We shall, however, follow the approach developed in~\cite{CFK}.
\begin{defn}
For $\lambda\in P^+$,   the {\em global Weyl module} $W_A(\lambda)$
is the left $\bu(\lie g\otimes A)$--module generated by an element
$w_\lambda$ with defining relations,$$(\lie n^+\otimes A)w_\lambda=
0,\qquad hw_\lambda= \lambda(h)w_\lambda,\qquad
(x_{\alpha_i}^-)^{\lambda(h_{\alpha_i})+1}w_\lambda=0,$$ where
$h\in\lie h$ and $i\in I$. \hfill\qedsymbol\end{defn}
It is clear that $W_A(\lambda)$ is not isomorphic to $W_A(\mu)$ if $\lambda\not=\mu$.

Suppose that  $\varphi$ is an algebra automorphism of $A$.
The defining ideal of~$W_A(\lambda)$ is clearly preserved by the automorphism of $\bu(\lie g\otimes A)$ induced by 
the Lie algebra automorphism $1\tensor\varphi$ of $\lie g\tensor A$ (cf.~\ref{subs:prelim.20}), and  we have an isomorphism of $\lie g\otimes A$--modules,
\begin{equation}\label{invaut}W_A(\lambda)\cong (1\tensor\varphi)^* W_A(\lambda).\end{equation}

\subsection{} The following construction shows immediately that $W_A(\lambda)$ is non--zero.
 Given any ideal $\mathfrak I$ of $A$, define an action of $\lie g\otimes A$ on
$V(\lambda)\otimes A/\mathfrak I$ by
$$(x\otimes a)(v\tensor b)=xv\otimes\bar a b,\qquad x\in \lie g,\, a\in A,\, b\in A/\mathfrak I,$$ where $\bar a$ is the canonical image of $a$ in $A/\mathfrak I$.
In particular, if $a\notin\mathfrak I$ and $h\in\lie h$ is such that $\lambda(h)\ne 0$  we have
\begin{equation}\label{nonzero}(h\otimes a)(v_{\lambda}\otimes 1)=\lambda(h)v_{\lambda}\otimes \bar
a\ne 0.\end{equation}
Clearly if $\mathfrak I\in\Max A$, then
$$V(\lambda)\otimes A/\mathfrak I\cong V(\lambda)$$ as $\lie
g$--modules and hence  $v_{\lambda}\otimes 1$ generates $V(\lambda)\tensor A/\mathfrak I$ as a $\lie g$--module (and so also as a $\lie g\otimes A$--module). Since $v_\lambda\otimes 1$ satisfies
the defining relations of $W_A(\lambda)$, we see that  $V(\lambda)\tensor A/\mathfrak I$ is a non--zero quotient of $W_A(\lambda)$.

\subsection{} Given a weight module $V$ of $\lie g\otimes A$, and a Lie subalgebra $\lie a$ of $\lie g\otimes A$,
set $$V_\mu^{\lie a}=V_\mu\cap V^{\lie a},\qquad \mu\in\lie h^*.$$

The following lemma is elementary. \begin{lem}\label{elem} For $\lambda\in P^+$ the module $W_A(\lambda)$ is a locally finite--dimensional $\lie g$--module and we have $$W_A(\lambda)=\bigoplus_{\eta\in Q^+
}W_A(\lambda)_{\lambda-\eta}$$ which in particular means that $\wt W_A(\lambda)\subset\lambda-Q^+.$  If
 $V$ is a $\lie g\otimes A$--module which is  locally finite--dimensional as a $\lie g$--module then we have an isomorphism of vector spaces,
\begin{equation*}\Hom_{\lie
g\otimes A}(W_A(\lambda),V)\cong V_\lambda^{\lie n^+\otimes A}.\tag*{\qedsymbol}
\end{equation*}
\end{lem}

\subsection{}  The weight spaces  $W_A(\lambda)_{\lambda-\eta}$ are not necessarily finite--dimensional, and to understand them, we proceed as follows.
It is easily checked that one can regard $W_A(\lambda)$ as a right module for $\bu(\lie h\otimes A)$ by setting
\[(uw_\lambda)(h\otimes a) = u(h\otimes a)w_\lambda,\qquad u\in\bu(\lie g\otimes A),\ \ h\in\lie h,\ \ a\in A.\]
Since $\bu(\lie h\tensor A)$ is commutative, the algebra $\ba_\lambda$ defined by
\[ \ba_\lambda=\bu(\lie h\tensor A)/\ann_{\bu(\lie h\tensor A)}w_\lambda\]
is a commutative associative algebra.
It follows that  $W_A(\lambda)$ is a $(\lie g\otimes A,
\ba_\lambda)$-bimodule and that for all $\eta\in Q^+$, the weight space
$W_A(\lambda)_{\lambda-\eta}$ is a right $\ba_\lambda$--module. Clearly
$$W_A(\lambda)_\lambda\cong_{\ba_\lambda}\ba_\lambda,$$ where we
regard $\ba_\lambda$ as a right $\ba_\lambda$-module through the right
multiplication.
The following  is immediate.
\begin{lem} For $\lambda\in  P^+$, $\eta\in Q^+$ the  subspaces $W_A(\lambda)^{\lie n^+}_{\lambda-\eta}$
and $W_A(\lambda)^{\lie n^+\otimes A}_{\lambda-\eta}$ are $\ba_\lambda$--submodules of $W_A(\lambda)$ and we have
\begin{equation*}
W_A(\lambda)_\lambda^{\lie n^+} =W_A(\lambda)_\lambda^{\lie n^+\otimes
 A}=W_A(\lambda)_\lambda.
 \end{equation*} Further, if  $\mu\in P^+$, the space $\Hom_{\lie g\otimes A}(W_A(\mu), W_A(\lambda))$
  has the natural structure of a right $\ba_\lambda$--module and
  \begin{equation*}\Hom_{\lie g\otimes A}(W_A(\mu), W_A(\lambda))\cong_{\ba_\lambda}  W_A(\lambda)_\mu^{\lie n^+\otimes A}.\tag*{\qedsymbol}
  \end{equation*}
\end{lem}

\subsection{}

   \begin{lem}\label{ninv} For  $\lambda\in P^+$, $\boa\in\ba_\lambda$ the assignment
   $w_\lambda\to w_\lambda\boa$ extends to a well--defined homomorphism $W_A(\lambda)\to W_A(\lambda)$ of  $\lie g\otimes A$--modules and
 we have $$\Hom_{\lie g\otimes A}(W_A(\lambda),W_A(\lambda))\cong_{\ba_
 \lambda}\ba_\lambda\cong_{\ba_\lambda} W_A(\lambda)_\lambda^{\lie n^+\otimes A}. $$
 \end{lem}
    \begin{pf} Since $W_A(\lambda)$ is a $(\lie g\otimes A,\ba_\lambda)$--bimodule, it follows that $w_\lambda\boa$ satisfies the defining relations of $W_A(\lambda)$ which yields the first statement of the Lemma. For the second,
     let $\pi:W_A(\lambda)\to W_A(\lambda)$ be a nonzero $\lie g\otimes A$--module map.
     Since $W_A(\lambda)_\lambda =\bu(\lie h\otimes A)w_\lambda$, there exists $u_\pi\in\bu(\lie h\otimes A)$ such that $\pi(w_\lambda)= u_\pi w_\lambda$.
     Since $\pi$ is non--zero, the image $\tilde u_\pi$ of $u_\pi$ in $\ba_\lambda$ is non--zero. Thus, we obtain a well-defined map 
     $\Hom_{\lie g\otimes A}(W_A(\lambda),W_A(\lambda))\to\ba_\lambda$ given by $\pi\mapsto \tilde u_\pi$. It is clear that the map is an isomorphism of right $\ba_\lambda$--modules and the lemma is proved.
    \end{pf}

\section{The main results}\label{sect:main res}

\subsection{} For  $\bos=(s_i)_{i\in I}\in\bz_+^I$, set
\begin{equation}\label{fundtensor}\ba_\bos=\bigotimes_{i\in
I}\ba_{\omega_i}^{\otimes s_i}, \qquad
W_A(\bos)=\bigotimes_{i\in I} W_A(\omega_i)^{\tensor s_i},
\qquad w_\bos=\bigotimes_{i\in I}w_{\omega_i}^{\otimes s_i}, \end{equation}
where all tensor products are taken in the same (fixed) order.
Given $k,\ell\in\bz_+$ let $\cal R_{k,\ell}$ be the algebra of  polynomials $\bc[t_1^{\pm1},\dots, t_k^{\pm1}, u_1,\dots,u_\ell]$  with the convention that if $k=0$ (respectively $\ell=0$), $\cal R_{0,\ell}$ (respectively $\cal R_{k,0}$) is just the ring of polynomials (respectively Laurent polynomials) in $\ell$ (respectively $k$) variables.
\subsection{} The main result of this
paper is the following.
\begin{thm}\label{thm1} Assume that
 $A=\cal R_{k,\ell}$ for some $k,\ell\in\bz_+$.
 For all $\bos\in\bz_+^n$ and $\mu\in P^+$, we have
\begin{equation} \label{homiso} \Hom_{\lie g\otimes
 A}(W_A(\mu),W_A(\bos))\cong_{\ba_\bos} W_A(\bos)^{\lie n^+\otimes A}_\mu=\Big(\bigotimes_{i\in I}(W_A(\omega_i)^{\lie n^+\otimes A})^{\otimes s_i}\Big)_\mu.\end{equation}
 \end{thm}

In the case when $\lie g$ is a classical
simple Lie algebra, we can make \eqref{homiso} more precise. Let
$I_0$ be the set of~$i\in I$ such that $\alpha_i$ occurs in~$\theta$ with the coefficient $2/(\alpha_i,\alpha_i)$. In particular,
$I_0=I$ for $\lie g$ of type~$A$ or~$C$.
Given
$\bos=(s_i)_{i\in I}\in\bz_+^I$ and $\lambda\in P^+$, let $\boc_\bos(\lambda)\in\bz_+$ be the coefficient of  
 $e(\lambda)$ in
$$ \prod_{i\in I_0}e(\omega_i)^{s_i}\prod_{i\notin
 I_0}\left(\sum_{0\le j\le i/2 }\binom{j+k-1}{j}e(\omega_{i-2j})\right)^{s_i},$$
 where $\omega_0=0$.
\begin{cor}\label{class}
Let $\lambda\in P^+$ and $\bos\in\bz_+^I$.
Assume either that $\lie g$ is not an exceptional Lie algebra or that $s_i=0$ if $i\notin I_0$.
 We have
 $$
 \Hom_{\lie g\otimes A}(W_A(\lambda),
 W_A(\bos))\cong_{\ba_\bos}\ba_\bos^{\oplus\boc_\bos(\lambda)},$$
where we use the convention that $\ba_\bos^{\oplus\boc_\bos(\lambda)}=0$ if $\boc_\bos(\lambda)=0$.
\end{cor}

\subsection{} Our next result is the following.  Recall from Lemma \ref{ninv} that for all $\lambda\in P^+$ we have $\Hom_{\lie g\otimes A}(W_A(\lambda), W_A(\lambda))\cong_{\ba_\lambda}\ba_\lambda$.
\begin{thm}\label{thm2}
 Let $A$ be the ring $\cal R_{0,1}$ or $\cal R_{1,0}$. For all $\mu=\sum_{i\in I}s_i\omega_i\in P^+$ with $s_i=0$ if $i\notin I_0$, and all $\lambda\in P^+$, we have
 $$
 \Hom_{\lie g\otimes A}(W_A(\lambda),W_A(\mu))=0, \qquad\text{if $\lambda\ne\mu$}.
 $$
Further, any non--zero element of $\Hom_{\lie g\otimes A}(W_A(\mu), W_A(\mu))$ is injective.
An analogous result holds when $A =\cal R_{k,\ell}$, $k,\ell\in\bz_+$,  $\lie g\cong\lie {sl}_{n+1}$ and
 $\mu=s\omega_1$.
  \end{thm}

\subsection{} We now make some comments on the various restrictions in the main results. The proof of Theorem \ref{thm1} relies on an explicit construction of the fundamental global Weyl modules in terms of certain finite--dimensional modules called the fundamental local Weyl modules. 
A crucial ingredient of this construction is a natural bialgebra structure of 
$\cal R_{k,\ell}$. 
The proof of Corollary \ref{class} depends on a deeper understanding of the $\lie g$--module structure of the local fundamental Weyl modules. These results are unavailable for the exceptional Lie algebras when $k+\ell>1$. In the case when $k+\ell=1$, the structure of these modules for the exceptional algebras is known as a consequence of the work of many authors on the Kirillov--Reshetikhin conjecture (see \cite{CH} for extensive references on the subject). Hence, a precise statement of Corollary~\ref{class}  could be made when $k+\ell=1$ in a case by case and in a not very compact fashion. The interested reader is referred to~\cite{HKOTY} and~\cite{Kl}.

\subsection{}  Before discussing Theorem \ref{thm2}, we make the following conjecture.
\begin{conj}
\label{injective2} Let $A=\cal R_{k,\ell}$ for some $k,\ell\in\bz_+$.
  Then for all $\lambda\in P^+$ and $\bos\in\bz_+^n$, any non--zero element of  $\Hom_{\lie g\otimes A}(W_A(\lambda), W_A(\bos))$ is injective.
\end{conj}

  The proof of Theorem~\ref{thm2} will rely on the fact that this conjecture is true (see Section 5  of this paper) when $k+\ell=1$ and $s_i=0$ if $i\notin I_0$
  {\em as well as} on the fact that the fundamental local Weyl module is irreducible as a $\lie g$--module if $i\in I_0$. We shall prove in Section 5  using Corollary \ref{class} and the work of
  \cite{FL} that the conjecture is also true when $\lie g=\lie{sl}_{r+1}$ and $\lambda=s\omega_1$.  Remark \ref{red} of this paper shows that  $\Hom_{\lie g\otimes A}
  (W_A(\lambda),W_A(\mu))$ can be non--zero if we remove the restriction on $\mu$.

\subsection{}  Finally, we make some remarks on quantum analogs of this result. In the case of the  quantum loop algebra, one also has analogous notions  of global and local Weyl modules which were  defined in  \cite{CPweyl}, and one can construct the global fundamental Weyl module from the local Weyl module in a way analogous to the one given in this paper for $A=\bc[t,t^{-1}]$. It was shown in \cite{CPweyl} for the quantum loop algebra of~$\lie{sl}_2$
that the canonical map from the global Weyl module into the tensor product of fundamental global Weyl modules is injective.
For the general quantum loop algebra, this was established by
Beck and Nakajima (\cite{BN}) using crystal and global bases. They also describe the space of extremal weight vectors in the tensor product of quantum fundamental global Weyl modules.

\section{The algebra \texorpdfstring{$\ba_\lambda$}{A\_lambda} and the local Weyl modules}\label{sect:local Weyl}  In this section we recall some necessary results from \cite{CFK} and also the definition and elementary properties of local Weyl modules.
\subsection{} For  $r\in\bz_+$
 the symmetric group $S_r$ acts naturally on $A^{\otimes r}$ and on $(\Max A)^{\times r}$ and we let
 $(A^{\otimes r})^{S_r}$ be the corresponding ring of invariants and $(\Max A)^{\times r}/S_r$ the set of orbits.
 If $r=r_1+\cdots +r_n$, then
  we regard $S_{r_1}\times \cdots\times S_{r_n}$ as a subgroup of $S_r$ in the canonical way, i.e. $S_{r_1}$ permutes the first $r_1$ letters, $S_{r_2}$ the
  next $r_2$ letters and so on. Given $\lambda=\sum_{i\in I}r_i\omega_i\in P^+$, set
  \begin{gather}\label{defbba} r_\lambda=\sum_{i\in I} r_i,\ \ S_\lambda=S_{r_1}\times\cdots\times S_{r_n},\ \
  \ \ \mathbb{A}_\lambda=(A^{\otimes r_\lambda})^{S_\lambda},\\ \Max \mathbb A_\lambda=(\Max A)^{r_\lambda}/S_\lambda.
 \end{gather}
 The algebra $\mathbb A_\lambda$ is generated by elements of the form  \begin{equation}\label{defnsym} \sym^i_\lambda(a)=1^{\otimes(r_1+\cdots+r_{i-1})}\otimes \Big(\sum_{k=0}^{r_i-1} 1^{\otimes k} \otimes a\otimes 1^{\otimes(r_i-k-1)}\Big)\otimes 1^{\otimes (r_{i+1}+\cdots+ r_n)}, \qquad a\in A,\, i\in I.\end{equation}

 The following was proved in \cite[Theorem~4]{CFK}.
\begin{prop}\label{bastru} Let $A$ be a finitely generated commutative associative algebra over $\bc$ with trivial Jacobson radical.  
Then the homomorphism of associative algebras $\bu(\lie h\otimes A)\to \mathbb A_\lambda$ defined by
\begin{equation*}
h_i\otimes a\mapsto \sym^i_\lambda(a),\ \ i\in I, \ \ a\in A
\end{equation*}
induces an isomorphism of algebras
$\sym_\lambda:\ba_\lambda\stackrel\sim\longrightarrow \mathbb
A_\lambda$. In particular, if $A$ is a finitely generated integral
domain then $\ba_\lambda$ is isomorphic to an integral subdomain of
$\ba_\bor$. \hfill\qedsymbol
\end{prop}
\subsection{}\label{ahom}  For  $\lambda,\mu\in P^+$, it is clear that the tensor
product $W_A(\lambda)\otimes W_A(\mu)$ has the natural structure of a $(\lie g\otimes A,\ba_\lambda\otimes\ba_\mu)$--module. We recall from  \cite{CFK} that, in fact, there exists a   $(\lie
g\otimes A,\ba_{\lambda+\mu})$--bimodule structure on $W_A(\lambda)\otimes W_A(\mu)$.

It is clear from Definition~\ref{defn-gl-weyl} that the assignment $w_{\lambda+\mu}\mapsto w_\lambda\tensor w_\mu$ defines a homomorphism
$\tau_{\lambda,\mu}:W_A(\lambda+\mu)\to
W_A(\lambda)\otimes W_A(\mu)$ of  $\lie g\otimes A$--modules.
The restriction of this map to $W_A(\lambda+\mu)_{\lambda+\mu}$ induces a homomorphism  of algebras
$\ba_{\lambda+\mu}\to\ba_\lambda\otimes\ba_\mu$ as follows. Consider
the restriction of the comultiplication $\Delta$  of $\bu(\lie g\otimes A)$
to  $\bu(\lie h\otimes A)$. It is not hard to see that
 \begin{multline*}
 \ann_{\bu(\lie h\tensor A)\tensor\bu(\lie h\tensor A)}(w_\lambda\tensor w_\mu)
 \subset\ann_{\bu(\lie h\otimes A)}w_\lambda\otimes\bu(\lie h\otimes A)+\bu(\lie h\otimes A)\otimes\ann_{\bu(\lie h\otimes
 A)}w_\mu,\end{multline*}
 and hence we have  $$\Delta(\ann_{\bu(\lie h\otimes A)}(w_{\lambda+\mu}))\subset\ann_{\bu(\lie h\otimes
 A)}w_\lambda\otimes\bu(\lie h\otimes A)+\bu(\lie h\otimes A)\otimes\ann_{\bu(\lie h\otimes
 A)}w_\mu.$$
 It is now immediate that the comultiplication $\Delta:\bu(\lie h\tensor A)\to \bu(\lie h\tensor A)\tensor \bu(\lie h\tensor A)$
induces a homomorphism of algebras
$\Delta_{\lambda,\mu}:
 \ba_{\lambda+\mu}\to\ba_\lambda\otimes\ba_\mu$. This endows any right $\ba_\lambda\tensor
 \ba_\mu$-module (hence, in particular, $W_A(\lambda)\tensor W_A(\mu)$) with the structure of a right $\ba_{\lambda+\mu}$-module.
It was shown in \cite{CFK} that
$\tau_{\lambda,\mu}$ is then a homomorphism of $(\lie g\tensor A,\ba_{\lambda+\mu})$-bimodules. Summarizing, we
have
\begin{lem} Let $\lambda_s\in P^+$, $1\le s\le k$ and let
$\lambda=\sum_{s=1}^k\lambda_s$. The natural map $W_A(\lambda)\to
W_A(\lambda_1)\otimes \cdots\otimes W_A(\lambda_k)$ given by
$w_\lambda\mapsto w_{\lambda_1}\otimes\cdots\otimes w_{\lambda_k}$ is
a homomorphism of $(\lie g\otimes A,\ba_\lambda)$--bimodules.\qed\end{lem}
\subsection{} For~$\lambda\in P^+$,
let~$\operatorname{mod} \ba_\lambda$ be the category of left $\ba_\lambda$--modules and let $\bw^\lambda_A $ be the
  right exact functor from~$\operatorname{mod}\ba_\lambda$ to
   the category of $\lie g\otimes A$--modules given
    on objects  by $$\bw_A^\lambda M=W_A(\lambda)\otimes_{\ba_\lambda} M,\qquad M\in\operatorname{mod} \ba_\lambda.$$
     Clearly $\bw^\lambda_AM$ is a weight module for $\lie g$  and we have isomorphisms of vector spaces \begin{gather*}(\bw^\lambda_AM)_{\lambda-\eta}\cong (W_A(\lambda))_{\lambda-\eta}\otimes_{\ba_\lambda}M,\qquad \eta\in Q^+,\\
(\bw^\lambda_AM)_{\lambda}\cong
(W_A(\lambda)_\lambda)\otimes_{\ba_\lambda}M\cong
w_\lambda\otimes_{\bc} M.\end{gather*} Moreover, $\bw^\lambda_AM$ is
generated as a $\lie g\otimes A$--module by the space
$w_\lambda\otimes_\bc M$ and
 $$\bw^\lambda_AM\cong\bw^\mu_AN\iff\ \ \lambda=\mu,\ \ M\cong_{\ba_\lambda}N.$$ The isomorphism classes of simple objects of~$\operatorname{mod} \ba_\lambda$ are given by the
 maximal ideals  of $\ba_\lambda$. Given~$\bi\in\Max\ba_\lambda$, the quotient
$\ba_\lambda/\bi$ is a simple object of $\operatorname{mod}\ba_\lambda$ and has
dimension one. The $\lie g\otimes A$--modules $\bw^\lambda_A\ba_\lambda/\bi$ are called the local Weyl modules and when $\lambda=\omega_i$, $i\in I$ we call them the fundamental local Weyl modules.
It follows that  $(\bw^\lambda_A\ba_\lambda/\bi)_\lambda$ is
also a  one--dimensional vector space spanned by $$w_{\lambda, \ba_\lambda/\bi}=w_\lambda\otimes 1.$$  We note the following corollary.
\begin{cor} Suppose that $M\in\operatorname{mod}\ba_\lambda$ is
finite--dimensional. Then $\bw^\lambda_AM$ is finite-dimen\-sional.
In particular, the local Weyl modules are finite--dimensional and
have a unique irreducible quotient $\bv^\lambda_AM$.
\end{cor}

\subsection{}\label{haction} We note the following consequence of Proposition \ref{bastru}.
\begin{lem} For $i\in I$, we have $\ba_{\omega_i}\cong A$  and  $W_A(\omega_i)$ is a  finitely
generated  right $A$--module.  The fundamental local Weyl modules are given by
\begin{equation*}
\bw^{\omega_i}_AA/\mathfrak I =W_A(\omega_i)\otimes_A A/\mathfrak
I,\qquad\mathfrak I\in\Max A.\end{equation*} In particular, we have \begin{equation*} (h\otimes a)w_{\omega_i, A/\mathfrak I}=0,\quad (h\otimes b)w_{\omega_i, A/\mathfrak I}=\omega_i(h)w_{\omega_i}\otimes \bar b,\qquad h\in \lie h,\, a\in\mathfrak I,\,b\in A. \tag*{\qedsymbol}\end{equation*}
\end{lem}

\subsection{}\label{lem.fingen}  The following lemma is a special case of a result proved in \cite{CFK}
 and we include the proof in this case for the reader's convenience.
\begin{lem} Let $A$ be a finitely generated, commutative associative algebra. For  $\mathfrak I\in\Max A$ there exists $N\in\bz_+$ such that for all $i\in I$,
\begin{equation}\label{fingen}(\lie g\otimes \mathfrak I^N)\bw_A^{\omega_i}A/\mathfrak I=0.\end{equation}
\end{lem}
\begin{pf} Since $\lie g$ is a simple Lie algebra
  \eqref{fingen} follows if we show that there exists $N\ge0$ such that
  $$(x^-_\theta\otimes \mathfrak I^N)\bw_A^{\omega_i}A/\mathfrak
I=0.$$
Since $$\bw_A^{\omega_i}A/\mathfrak
I=\bu(\lie n^-\otimes A)w_{\omega_i,A/\mathfrak I},\qquad
   [x^-_\theta,\lie n^-]=0,$$  it is enough
   to show that there exists $N\ge0$ such that
    $$(x^-_\theta\otimes\mathfrak I^N)w_{\omega_i,A/\mathfrak I}=0.$$
     We claim that is a consequence of showing that for $j\in I$
     there exists $N_j\in\bz_+$, with
      \begin{equation}\label{simplecase}(x^-_{\alpha_j}\otimes \mathfrak I^{N_j})w_{\omega_i,A/\mathfrak I}=0.\end{equation}
     To prove the claim, write
     $x^-_\theta=[x^-_{\alpha_{i_1}}[\cdots[x^-_{\alpha_{i_{p-1}}},x^-_{\alpha_{i_p}}]\cdots] $
     for some $i_1,\dots, i_p\in I$ and take $N=\sum_{j=1}^pN_{i_j}
       $. To prove \eqref{simplecase}, observe first that for~$j\not=i\in I$, $k\in I$ and for all $a\in A$
$$
x^+_{\alpha_k}(x^-_{\alpha_j}\tensor a)w_{\omega_i,A/\mathfrak I}=\delta_{k,j}(h_j\tensor a)w_{\omega_i,\mathfrak I}=0,
$$
by Lemma~\ref{haction} and the defining relations of~$W_A(\omega_i)$. Thus, $(x^-_{\alpha_j}\tensor a)w_{\omega_i,A/\mathfrak I}\in
(\bw_A^{\omega_i}A/\mathfrak I)^{\lie n^+}$ and since $\omega_i-\alpha_j\notin P^+$ we conclude that
        $$(x^-_{\alpha_j}\otimes A)w_{\omega_i,A/\mathfrak I}=0,\qquad j\ne i.$$
If $j=i$, then $$0=(x_{\alpha_i}^+\otimes
a)(x^-_{\alpha_i})^2
         w_{\omega_i}=2((x^-_{\alpha_i}\otimes 1)(h_i\otimes a)-(x^-_{\alpha_i}\otimes a))
         w_{\omega_i}.$$  By Lemma~\ref{haction}, we have $(h_i\otimes a)w_{\omega_i,A/\mathfrak
         I}=0$ if $a\in\mathfrak I$, and so  we get
          $$(x^-_{\alpha_i}\otimes a)w_{\omega_i,A/\mathfrak I}=0,\qquad a\in\mathfrak
          I$$ and
          \eqref{simplecase} is established.
\end{pf}

\section{Fundamental global Weyl modules}\label{sect:Loc and glob Weyl} In this section we establish the main tool for proving Theorem \ref{thm1}.   It is not, in general, clear how   (or even if it is possible)
to reconstruct the global Weyl module from a local Weyl module. The main result of this section is that it is possible to do so when $\lambda=\omega_i$ and $A=\cal R_{k,\ell}$ for some $k,\ell\in\bz_+$.

\subsection{}\label{subs:gen-coalg} We begin with a general construction. The Lie algebra $(\lie g\otimes A)\otimes A$ acts naturally on $V\otimes A$ for any $\lie g\otimes A$--module $V$.
 Suppose  that  $A$ is a  bialgebra  with the comultiplication $\boh:A\to A\otimes
 A$. (It is useful to recall that $A$ is a commutative associative algebra with identity). Then the  comultiplication map $\boh$
 induces a homomorphism of Lie algebras $1\tensor\boh: \lie g\tensor A\to \lie g\tensor A\otimes A$ (cf.~\ref{subs:prelim.20})
 and thus a $\lie g\tensor A$-module structure on~$V\tensor A$. Explicitly, the $(\lie g\otimes A,A)$-bimodule
  structure on $V\otimes A$ is given by the following formulas:
  $$(x\otimes a)(v\otimes b)=\sum_{s}(x\otimes a_s')v\otimes a_s''b,\qquad (v\otimes b)a=v\otimes ba, \, v\in V, a,b\in A,$$
   where $\boh(a)=\sum_sa_s'\otimes a_s''$.  We denote this  bimodule by $(V\otimes A)_\boh$ and observe that it is a  free right $A$--module of rank equal to $\dim_\bc V$.
   It is trivial to see that  $(V\otimes A)_\boh$ is a weight module for $\lie g\otimes A$ if $V$ is a weight module for $\lie g\otimes A$ and that $$((V\otimes A)_\boh)_\mu= V_\mu\otimes A,\qquad  V^{\lie n^+\otimes A}\otimes A\subset (V\otimes A)_\boh^{\lie n^+\otimes A}.$$
   Moreover, if $V_1,V_2$ are $\lie g\otimes A$--modules, one has a natural inclusion
   \begin{equation} \label{inc} \Hom_{\lie g\otimes A}(V_1,V_2)\hookrightarrow\Hom_{\lie g\otimes A}((V_1\otimes A)_\boh, (V_2\otimes A)_\boh),\ \ \eta\mapsto \eta\otimes 1.\end{equation} In particular, if $V$ is reducible, then $(V\otimes A)_\boh$ is also a reducible $(\lie g\otimes A)$--module.

\subsection{}\label{lem:comult}
Let $\boh_{k,\ell}:\cal R_{k,\ell}\to \cal R_{k,\ell}\otimes \cal R_{k,\ell}$ be the  comultiplication given by,
$$
 \boh_{k,\ell}(t_s^{\pm 1})=t_s^{\pm 1}\tensor t_s^{\pm 1},\qquad \boh_{k,\ell}(u_r)=u_r\tensor 1+1\tensor u_r,\qquad
$$
where $1\le s\le k$ and  $1\le r\le \ell$.
Any monomial $\bom\in \mathcal R_{k,\ell}$ can be written uniquely as a product of monomials
$$\bom=\bom_u \bom_t,\qquad\bom_t\in\bc[t_1^{\pm 1},\dots, t_k^{\pm 1}],\,
\bom_u\in\bc[u_1,\dots, u_\ell].$$  Set $\deg t_s^{\pm 1}= \pm1$ and $\deg u_r=1$ for $1\le s\le k$, $1\le r\le\ell$ and
let  $\deg_t \bom$ (respectively, $\deg_u \bom$)  be the total degree of~$\bom_t$ (respectively, $\bom_u$) and for  for any $f\in A$ define  $\deg_tf$ and $\deg_u f$ in the obvious way.
The next lemma is elementary.
\begin{lem}\label{monomials} Let $\bom=\bom_t\bom_u$ be a monomial in $\cal R_{k,\ell}$. Then $\bom_t\in \mathcal R_{k,\ell}^\times$, $\boh_{k,\ell}(\bom_t)=
\bom_t\tensor \bom_t$ and
\begin{equation}\label{eq:boh-mon}
\boh_{k,\ell}(\bom)=\bom\tensor \bom_t+\sum_{q} \bom'_{u,q}\bom_t\tensor \bom''_{u,q} \bom_t=
\bom_t\tensor \bom+\sum_q \bom''_{u,q}\bom_t\tensor \bom'_{u,q}\bom_t,
\end{equation}
where  $\bom'_{u,q},\bom''_{u,q}$ are (scalar multiples of) monomials in the $u_r$, $1\le r\le \ell$, such that if $\bom'_{u,q}\ne 0$,
$\bom''_{u,q}\ne 0$, then $\deg_u \bom_{u,q}'<\deg_u\bom$ and $\deg_u\bom_{u,q}'+\deg_u\bom_{u,q}''=\deg_u \bom$.\qed
\end{lem}

\subsection{}  {\em For the rest of the section $A$ denotes the algebra $\cal R_{k,\ell}$ for some $k,\ell\in\bz_+$ and $\mathfrak I$ the ideal of $A$ generated by the elements $\{t_1-1,\dots, t_k-1, u_1,\dots, u_\ell\}$.}

Suppose that    $\mathfrak J\in\Max \cal R_{k,\ell}$. It is clear that there exists an algebra automorphism $\varphi:\cal R_{k,\ell}\to \cal R_{k,\ell}$ such that $\varphi(\mathfrak I)=\mathfrak J$.  As a consequence, we have an induced isomorphism  of $\lie g\otimes A$--modules, \begin{equation}\label{indep0}\bw^{\omega_i}_AA/\mathfrak I\cong (1\otimes \varphi)^*\bw^{\omega_i}_A A/\mathfrak J. \end{equation}Moreover, if we set
$$\boh_{k,\ell}^\varphi =(\varphi\tensor\varphi)\circ\boh_{k,\ell}\circ\varphi^{-1}:A\to A\tensor A,$$  then $\boh_{k,\ell}^\varphi$ also defines a
bialgebra structure on $A$ and we have an isomorphism of $\lie g\otimes A$-modules  \begin{equation}\label{indep} (1\otimes \varphi)^*(\bw^{\omega_i}_AA/\mathfrak I\otimes A)_{\boh_{k,\ell}}\cong (\bw^{\omega_i}_AA/\mathfrak J\otimes A)_{\boh_{k,\ell}^\varphi}.\end{equation}
This becomes an isomorphism of $(\lie g\tensor A,A)$-bimodules if we twist the right $A$-module structure of $(\bw^{\omega_i}_AA/\mathfrak J\otimes A)_{\boh_{k,\ell}^\varphi}$
by~$\varphi$.

\subsection{}\label{thm:glfundreal} We now reconstruct the global fundamental Weyl module from a local one. \begin{prop}\label{glfundreal}
 For all $i\in I$,
  the assignment $w_{\omega_i}\mapsto w_{\omega_i,A/\mathfrak I}\tensor 1$
  defines an  isomorphism of $(\lie g\otimes A,A)$--bimodules
   $$W_A(\omega_i)\stackrel{\cong}{\longrightarrow}(\bw^{\omega_i}_AA/\mathfrak I\otimes A)_{\boh_{k,\ell}}.$$
\end{prop}
\begin{rem}  It is clear from \eqref{indep} and Section~\ref{defn-gl-weyl} that one can work with an arbitrary ideal $\mathfrak J$ provided that $\boh_{k,\ell}$ is replaced by $\boh_{k,\ell}^\varphi$, where $\varphi$ is the unique automorphism of~$A$ such that $\varphi(\mathfrak I)=\mathfrak J$. \end{rem}
\begin{pf}
The element $w_{\omega_i,A/\mathfrak I}\tensor 1\in (\bw_A^{\omega_i} A/\mathfrak I\tensor A)_\boh$
satisfies the relations in Definition~\ref{defn-gl-weyl} and  hence the assignment
$w_{\omega_i}\mapsto w_{\omega_i,A/\mathfrak I}\tensor 1$ defines a homomorphism of $\lie g\tensor A$-modules
$$\bop: W_A(\omega_i)\to (\bw_A^{\omega_i} A/\mathfrak I\tensor A)_{\boh_{k,\ell}}.
$$
We begin by proving that $\bop$ is a homomorphism of right $A$-modules.
 Using Lemma~\ref{haction} and the definition of the right module structure on $W_A(\omega_i)$ we see that
\begin{equation*}
\bop((uw_{\omega_i})a)=\bop(u(h_i\otimes
   a)w_{\omega_i})= u(h_i\otimes a)(w_{\omega_i,A/\mathfrak I}\otimes 1),\qquad u\in\bu(\lie g\otimes A), \, a\in A.
\end{equation*}
This shows that  it is enough to prove that for any monomial~$\bom$ in~$A$, we have
    \begin{equation}\label{bimap2}(h_i\otimes \bom)(w_{\omega_i,A/\mathfrak I}\otimes
   1)=
   w_{\omega_i,A/\mathfrak I}\otimes \bom. \end{equation}
Write $\bom=\bom_t\bom_u$ and observe that $\bom_t-1\in \mathfrak I$
while
$$
\deg_u\bom>0\implies\bom\in\mathfrak I.
$$
Using~\eqref{eq:boh-mon} we get
$$
\boh_{k,\ell}(\bom)-1\tensor \bom\in\mathfrak I\tensor A
$$
and since $(h_i\otimes \mathfrak I\otimes A)(w_{\omega_i,A/\mathfrak I}\otimes 1)=0$ by Lemma~\ref{haction},
we have
established~\eqref{bimap2}.

To prove that $\bop$ is surjective we must show that  $$\bu(\lie g\otimes A)(w_{\omega_i,A/\mathfrak I}\otimes 1)=(\bw^{\omega_i}_AA/\mathfrak I\otimes A)_{\boh_{k,\ell}}, $$ and   the remarks in  Section~\ref{subs:gen-coalg} show that it is enough to prove \begin{equation}\label{eq:tmp.2}
     (\bw^{\omega_i}_AA/\mathfrak I)_{\omega_i-\eta}\otimes A\subset \bu(\lie g\otimes A)(w_{\omega_i,A/\mathfrak I}\otimes 1),\qquad \eta\in Q^+.
\end{equation}
The argument is by induction on~$\Ht\eta$, the induction base with $\Ht\eta=0$ being immediate from~\eqref{bimap2}.
For the inductive step assume that we have proved the result for all $\eta\in Q^+$ with $\Ht\eta<k$. To prove the result for $\Ht\eta=k$, it suffices to prove that
for all $j\in I$ and  all monomials $\bom$ in $A$ we have
 $$((x^-_{\alpha_j}\otimes \bom)w)\otimes g\in \bu(\lie g\otimes A)(w_{\omega_i,A/\mathfrak I}\otimes 1),$$ where $w\in  (\bw^{\omega_i}_A A/\mathfrak I)_{\omega_i-\eta'}$ with $\Ht\eta'=k-1$ and $g\in A$.
 For this,  we argue by a further  induction on~$\deg_u \bom$. If $\deg_u \bom=0$ then  $\bom=\bom_t$
and we have
$$
((x^-_{\alpha_j}\tensor \bom)(w)\tensor g=
(x^-_{\alpha_j}\tensor \bom)(w\tensor \bom^{-1}g)\in \bu(\lie g\otimes A)(\bw^{\omega_i}_A A/\mathfrak I)_{\omega_i-\eta'}.
$$
This proves that the induction on $\deg_u\bom$ starts.
If $\deg_u \bom>0$ we  use~\eqref{eq:boh-mon} to get
$$
((x^-_{\alpha_j}\tensor \bom)w)\tensor g=(x^-_{\alpha_j}\tensor \bom)(w\tensor \bom_t^{-1} g)-\sum_q ((x^-_{\alpha_j}\tensor \bom_{u,q}'\bom_t)w)\tensor \bom_{u,q}'' g.
$$
Both terms on the right hand side are in $\bu(\lie g\otimes A)(\bw^{\omega_i}_A A/\mathfrak I)_{\omega_i-\eta'}$: the first by the induction hypothesis on $\Ht\eta'$
and the second by the induction hypothesis on~$\deg_u \bom$. This completes the proof of the surjectivity of~$\bop$.

To prove that $\bop$ is injective, recall from Section~\ref{subs:gen-coalg} that $(\bw^{\omega_i}_AA/\mathfrak I\otimes A)_{\boh_{k,\ell}}$
is a  free right $A$--module  of rank equal to the dimension of $\bw^{\omega_i}_AA/\mathfrak I$.
 Hence if $K$ is  the kernel of~$\bop$ we have
an isomorphism of right $A$-modules,
$$W_A(\omega_i)\cong K\oplus (\bw^{\omega_i}_AA/\mathfrak I\otimes A)_{\boh_{k,\ell}}.$$
Using \eqref{indep0} we see that for any maximal ideal~$\mathfrak J$ in~$A$,
$$
\dim (W_A(\omega_i)\tensor_A A/\mathfrak J)=\dim \bw^{\omega_i}_A A/\mathfrak J=
\dim\bw^{\omega_i}_A A/\mathfrak I\\=\dim((\bw^{\omega_i}_AA/\mathfrak I\otimes A)_{\boh_{k,\ell}} \tensor_A A/\mathfrak J).
$$
Therefore,
$K\tensor_A A/\mathfrak J=0$  and so by Nakayama's Lemma $K_{\mathfrak J}=0$
for all $\mathfrak J\in\Max A$. Thus, $K=0$.
\end{pf}

\section{Proof of Theorem \ref{thm1} and Corollary \ref{class}}\label{sect:pf Thm1}
  We continue to assume that $A=\cal R_{k,\ell}$, $k,\ell\in\bz_+$ and that $\mathfrak I$ is the  maximal ideal of $A$ generated by $\left\{t_1-1,\ldots, t_k-1,u_1,\ldots,u_\ell\right\}$. We also use the comultiplication $\boh_{k,\ell}$ and denote it by just $\boh$. Let
  $\mathfrak M_A\subset A$ be the set of monomials in the generators $u_r$, $t_s^{\pm 1}$, $1\le r\le \ell$, $1\le s\le k$.
 \subsection{} The following proposition, together with Lemma~\ref{lem.fingen}
 and Proposition~\ref{glfundreal}, completes the proof of Theorem~\ref{thm1}.

 \begin{prop}\label{invtensor}
 Suppose that  $V_s$, $1\le s\le M$ are  $\lie g\otimes A$--modules
 such that  there exists  $N\in\bz_+$ with \begin{equation}\label{defN}(\lie n^+\otimes\mathfrak I^N)V_s=0,\qquad 1\le s\le M.
 \end{equation}
Then
$$((V_1\otimes A)_\boh\otimes \cdots\otimes (V_M\otimes A)_\boh)^{\lie n^+\otimes A}=(V_1^{\lie n^+\otimes A}
\otimes A)_\boh\otimes\cdots\otimes (V_M^{\lie n^+\otimes A}\otimes
A)_\boh.$$
\end{prop}

\subsection{}\label{lem:lemI} The first step in the proof of Proposition~\ref{invtensor} is the following.
We need some notation. Let $V$ be a $\lie g\otimes A$-module and let $K\in\bz_+$. Define $$
V_{\ge K}=\{v\in V\,:\, (\lie n^+\otimes \bom) v=0,\, \bom\in\mathfrak M_A,\, |\deg_t\bom|\ge K\}.
$$ Note that $V_{\ge 0}= V^{\lie n^+\otimes A}.$
\begin{lem} Let $V$ be a $\lie g\otimes A$--module and  $K\in\bz_+$.
Then \begin{equation}\label{VK}((V\otimes A)_\boh)_{\ge K} = V_{\ge K}\otimes A.\end{equation}
In particular,
\begin{equation}\label{M1}(V\otimes A)_\boh^{\lie n^+\otimes A}=V^{\lie n^+\otimes A}\otimes A.\end{equation}
\end{lem}
\begin{pf}
Let $v_\boh\in (V\otimes A)_\boh$ and write, $v_\boh=\sum_p v_p\otimes g_p$, where $\{g_p\}_p$ is a linearly independent subset of $A$.
By~\eqref{eq:boh-mon} we have
$$
(x\otimes \bom)v_{\boh}=\sum_p (x\otimes \bom)v_p\otimes \bom_t g_p + \sum_{p,q} (x\tensor \bom'_{u,q} \bom_t )v_p\tensor \bom_{u,q}'' \bom_t g_p.
$$
with $\deg_u\bom'_{u,q}<\deg_u\bom$. Since $\deg_t\bom'_{u,q} \bom_t =\deg_t\bom$, it follows that $$V_{\ge K}\otimes A\subset ((V\otimes A)_\boh)_{\ge K}.$$
We prove the reverse inclusion  by induction on $\deg_u\bom$. Let $v_\boh\in ((V\tensor A)_\boh)_{\ge K}$ and let $\deg_t\bom\ge K$.
If $\deg_u \bom =0$,  then
$$0=(x\otimes \bom)v_{\boh}=\sum_p (x\otimes \bom)v_p\otimes g_p\bom .$$
Since the set  $\{g_p \bom\}_p$ is also linearly independent, we see that $(x\otimes \bom)v_p=0$ for all~$p$.
If $\deg_u \bom>0$, we use ~\eqref{eq:boh-mon} to get
$$
0=(x\otimes \bom)v_{\boh}=\sum_p (x\otimes \bom)v_p\otimes \bom_t g_p + \sum_{p,q} (x\tensor \bom'_{u,q} \bom_t )v_p\tensor \bom_{u,q}'' \bom_t g_p.
$$
Since $\deg_u\bom'_{u,q}<\deg_u\bom$ all terms in the second sum are zero by the induction hypothesis, and  the linear independence of the set $\{ \bom_t g_p\}_p$ gives $(x\otimes \bom)v_p=0$ for all~$p$.
\end{pf}

\subsection{}
\begin{prop} \label{crucial}
Let $U,V$ be $\lie g\tensor A$-modules  and suppose that for some $N\in\bz_+$
$$
(\lie n^+\tensor \mathfrak I^N)V=0.
$$
Then
\begin{equation}\label{M2case} (U\otimes (V\otimes A)_\boh)^{\lie n^+\otimes A}=
U^{\lie n^+\otimes A}\tensor (V^{\lie n^+\tensor A}\tensor A).\end{equation}
\end{prop}

Before proving this proposition, we establish Proposition~\ref{invtensor}. The argument is by induction on~$M$,
with \eqref{M1} showing that induction begins at $M=1$.  For~$M>1$,
take $$U=(V_1\otimes A)_\boh\otimes \cdots\otimes (V_{M-1}\otimes A)_\boh,\qquad
V= V_M. $$ The induction hypothesis applies to $U$  and together with Propsition \ref{crucial} completes the inductive step.

\subsection{}
\begin{lem}\label{lem:lemII} Let $A=\mathcal R_{k,\ell}$ with~$k>0$. Let $V$ be a $\lie g\tensor A$-module and
suppose that $(\lie n^+\otimes \mathfrak I^N)V=0$ for some $N\in\bz_+$.
Then for all $K\in\bz_+$
we have
$$
V^{\lie n^+\otimes A}=V_{\ge K}.
$$
\end{lem}
\begin{pf}
It suffices to prove that $V_{\ge K}\subset V_{\ge K-1}$ for all $K\ge 1$.
Since $(1-t_1^{\pm 1})^N\in\mathfrak I^N$ we have
\begin{equation}\label{reduction}0= (x\otimes\bom(1-t_1^{\pm 1})^N)v=(x\otimes \bom)v+\sum_{s=1}^N (-1)^s \binom Ns(x\tensor \bom t_1^{\pm s})v,\end{equation}for all $x\in\lie n^+$, $\bom\in \mathfrak M_A$ and $v\in V$
Suppose that $v\in V_{\ge K}$ and take $\bom\in\mathfrak M_A$ with $|\deg_t\bom|=K-1$. If~$\deg_t\bom\ge 0$ (respectively,
$\deg_t\bom<0$) then
$|\deg_t\bom t_1^s|\ge K$ (respectively, $|\deg_t\bom t_1^{-s}|\ge K$) for all~$s>0$.
Thus we conclude that all terms in the sum in~\eqref{reduction} with the appropriate sign choice equal zero hence
$(x\tensor \bom)v=0$ and so $v\in V_{\ge K-1}$.
\end{pf}

\subsection{}
\begin{lem}\label{lem:lemIII}
Let $A=\mathcal R_{0,\ell}$. Let $V$ be a $\lie g\tensor A$-module and
suppose that  $(\lie n^+\otimes \mathfrak I^N)V=0$ for some $N\in\bz_+$. Let  $K\ge N\in\bz_+$. Then
 \begin{equation}\label{biggenough2}
  V^{\lie n^+\otimes A}\otimes A=\{v_\boh\in (V\otimes A)_{\boh}\,:\, (\lie n^+\otimes \bom) v_\boh=0,\, \bom\in\mathfrak M_A,\,
  \deg_u\bom\ge K\}
 \end{equation}
\end{lem}
\begin{pf}
Since $(V\otimes A)_\boh$ is a $(\lie g\otimes A, A)$--bimodule the sets on both sides of \eqref{biggenough2} are right $A$--modules. Hence if $v_\boh$ is an element of the set on the right hand side of \eqref{biggenough2} then $v_\boh u_j^{s}$ is also in the right hand side of \eqref{biggenough2} for all $s\in\bz_+$.
Write $v_\boh=\sum_p v_p\tensor g_p$,
where $\{g_p\}_p$ is a linearly independent subset of~$A$.
Since the $u_j$, $1\le j\le \ell$ are primitive and $u_j^s\in\mathfrak I^N$ if~$s\ge N$, we have for all $0\le r\le N$
\begin{gather*} 0=(x\otimes
u_j^{(K+N-r)})(v_\boh)u_j^{(r)}=\sum_{s=0}^{N}\Big(\sum_p((x\otimes
u_j^{(s)})v_p)\otimes
u_j^{(K+N-r-s)}g_pu_j^{(r)}\Big)\\=\sum_{s=0}^N\binom{K+N-s}{r}\Big(\sum_p((x\otimes
u_j^{(s)})v_p)\otimes u_j^{(K+N-s)}g_p\Big).
\end{gather*}
We claim that the matrix $C(N,K)=(\binom{K+N-s}{r})_{0\le
s,r\le N}$ is invertible. Assuming the claim, we get
$$\sum_p((x\otimes
u_j^{(s)})v_p)\otimes u_j^{(K+N-s)}g_p=0,\qquad 0\le s\le N,$$ and  since the
$g_p$ are linearly independent this implies that $$(x\otimes
u_j^{(s)})v_p=0,\quad 0\le s\le N $$
and so $(x\tensor u_j^s)v_\boh=0$ for all $x\in\lie n^+$, $s\in\bz_+$.

Now, let $\bom\in\mathfrak M_A$ and let~$\alpha\in\Phi^+$. Then $(h_\alpha\tensor \bom)v_{\boh}$ is also
an element of the right hand side of~\eqref{biggenough2}
and hence by the preceding argument, we get \begin{align*}0&=
 (x_\alpha\otimes 1)(h_\alpha\otimes \bom)v_\boh\\&=
(h_\alpha\otimes \bom)(x_\alpha\otimes
1)v_\boh -2(x_\alpha\otimes \bom)v_\boh=
-(2x_\alpha\otimes \bom)v_\boh,\end{align*}  thus proving that $v_\boh\in (V\tensor A)_\boh^{\lie n^+\tensor A}=V^{\lie n^+\tensor A}\tensor A
$ by \eqref{M1}.

To prove the claim,
let $u$ be an indeterminate and let $\left\{p_r\in\C[u]: 0\le r\le N\right\}$ be a collection
of polynomials such that $\deg p_r=r$ (in particular, we assume that $p_0$ is a non-zero constant
polynomial). Then for any tuple $(a_0,\dots,a_N)\in \C^{N+1}$, we have
$\det(p_r(a_s))_{0\le r,s\le N}=c \det(a_s^r)_{0\le r,s\le N}=c\prod_{0\le r<s\le N} (a_s-a_r)$, where $c$ is the product of highest coeffcients
of the~$p_r$, $0\le r\le N$. Since $\binom{u}{r}$ is a polynomial in~$u$ of degree~$r$ with highest coefficient $1/r!$, we
obtain with $a_s=N+K-s$,
\begin{equation*}
\det C(N,K)=\Big(\prod_{r=1}^N r!\Big)^{-1} \prod_{0\le r<s\le N} (r-s)=(-1)^{N(N+1)/2}.\qedhere
\end{equation*}
\end{pf}

\subsection{}
Now we have all the necessary ingredients to prove Proposition~\ref{crucial}.
\begin{pf}[Proof of Proposition~\ref{crucial}]
Let $v_\boh\in (U\tensor (V\tensor A)_\boh)^{\lie n^+\tensor A}$ and write
$v_\boh=\sum_{p,s} w_s\otimes v_{s,p}\otimes g_p,$
where $\{w_s\}_s$ and $\{g_p\}_p$
are linearly independent subsets of $U$ and $A$ respectively.
 We have
\begin{align}\label{eq:tmp}
0&=(x\tensor \bom)v_\boh=\sum_{s,p} \Big(((x\tensor \bom)w_s)\tensor v_{s,p}\tensor g_p+w_s\tensor (x\tensor \bom)(v_{s,p}\tensor g_p)\Big)\\ &= \sum_{s,p}((x\tensor \bom)w_s)\tensor v_{s,p}\tensor g_p\nonumber\\&\quad+\sum_{s,p} w_s\otimes \Big( (x\otimes \bom)v_{s,p}\otimes g_p\bom_t+ \sum_q(x\otimes \bom'_{u,q}\bom_t)v_{s,p}\otimes \bom''_{u,q}\bom_tg_p\Big),\label{eq:tmp.a}
\end{align}
Suppose first that $A=\mathcal R_{k,\ell}$ with~$k>0$ and let $K=\max_p |\deg_t g_p|+1$.
If $\bom$ is such that $|\deg_t\bom|\ge K$, then the set $\{g_p\}_p$ is linearly independent from the set $\{\bom''_{u,q}\bom_tg_p\}_{p,q}$ and hence we must have that
\begin{equation}\label{eq:tmp.1}
\sum_{s,p} ((x\tensor \bom)w_s)\tensor v_{s,p}\tensor g_p=0,\qquad \sum_{s,p} w_s\tensor (x\tensor \bom)(v_{s,p}\tensor g_p)=0
\end{equation}
and using the linear independence of the elements $\{w_s\}_s$ we conclude that for all~$s$
$$
\sum_p (v_{s,p}\tensor g_p)\in ((V\tensor A)_\boh)_{\ge K}=V_{\ge K}\tensor A=V^{\lie n^+\tensor A}\tensor A,
$$
using~\eqref{VK} and Lemma~\ref{lem:lemII}. This proves Proposition~\ref{crucial} in the case when~$k>0$.

Suppose now that $A=\mathcal R_{0,\ell}$ and let $N\in\bz_+$ be such that $(\lie n^+\tensor \mathfrak I^{N})V=0$. Let $K=N+1+\max_p\left\{ \deg_u g_p\right\}$ and let $x\in\lie n^+$.
If~$\deg_u\bom\ge K$ then $\bom\in\mathfrak I^N$ and so
$(x\tensor\bom)v_{s,p}=0$. Furthermore,
$(x\tensor \bom'_{u,q})v_{s,p}\not=0$ implies that $\deg_u\bom'_{u,q}<N$. By Lemma~\ref{monomials} it follows that $\deg_u \bom''_{u,q}>\max_p\left\{ \deg g_p\right\}$. Therefore, the non-zero terms, if any, in the second
sum in~\eqref{eq:tmp.a} are linearly independent from those in the first sum and we obtain~\eqref{eq:tmp.1}. Furthermore,
we have
$$
0=\sum_{p}\sum_{\left\{q\,:\, \deg_u\bom'_{u,q}<N\right\}} (x\tensor \bom'_{u,q})v_{s,p}\tensor \bom''_{u,q} g_p
$$
and as before we conclude that $(x\tensor \bom'_{u,q})v_{s,p}=0$ when $\deg_u \bom'_{u,q}<N$. Thus, $(x\tensor\bom)\big(\sum_p v_{s,p}\tensor g_p\big)=0$
for all~$s$ and it remains to apply Lemma~\ref{lem:lemIII}.
\end{pf}

\subsection{} We conclude this section with a proof of Corollary \ref{class}.
The following is a special case of Theorem 7.1 and Proposition 7.7 of \cite{CFK}.
\begin{thm} \label{fundg} Let $\mathfrak I\in\Max A$, where $A=\mathcal R_{k,\ell}$. Then
\begin{gather*}\dim(\bw^{\omega_i}_AA/\mathfrak I)^{\lie n^+\otimes A}= \dim(\bw^{\omega_i}_AA/\mathfrak I)_{\omega_i}=1,\qquad i\in I_0,\end{gather*}
If  $\lie g$ is of type  $B_n$ or $D_n$ and $i\notin I_0$ then
\begin{gather*}
\dim(\bw^{\omega_i}_AA/\mathfrak I)^{\lie n^+\otimes A}_{\mu}=\begin{cases} 0,&\mu\ne \omega_{i-2j},\\
\displaystyle\binom{j+k-1}{j},&\mu=\omega_{i-2j},\,i-2j\ge 0\end{cases}\end{gather*} where $\omega_0=0$.
\qed\end{thm}

Another way to formulate this result is the following. The subspace $ (\bw^{\omega_i}_AA/\mathfrak I)^{\lie n^+\otimes A}$ is an $\lie h$--module with character given by
$$\ch (\bw^{\omega_i}_AA/\mathfrak I)^{\lie n^+\otimes A}
=\sum_{j: i-2j\ge 0}\binom{j+k-1}{j}e(\omega_{i-2j}).$$  Hence, using  Theorem \ref{thm1}, we get $$\ch  \bw_A(\bos)^{\lie n^+\otimes A}=\prod_{i\in I}(\ch (\bw^{\omega_i}_AA/\mathfrak I)^{\lie n^+\otimes A})^{s_i},$$
and  Corollary \ref{class} follows.

\section{Proof of Theorem \ref{thm2}}\label{sect:pf thm2}

Given $\bos\in\bz_+^I$, define $\mu_\bos\in P^+$ by
$\mu_\bos=\sum_{i\in I}s_i\omega_i$ and let $\tau_\bos:W_A(\mu_\bos)\to W_A(\bos)$
be the natural map of $(\lie g\otimes A,A)$--bimodules defined in the above lemma and
satisfying $\tau_\bos(w_{\mu_\bos})=w_\bos$.
Since $W_A(\bos)_{\mu_\bos}=w_\bos\otimes\ba_\bos$,
we see that any non--zero element $\tau$ of $\Hom_{\lie g\otimes A}(W_A(\mu_\bos),W_A(\bos))$
 is given by composing $\tau_\bos$  with right multiplication by an element of $\ba_\bos$,
  i.e. $\tau=\tau_\bos\boa$ with $\boa\in\ba_\bos$.

\subsection{}
\begin{lem}\label{conseq1} Assume that   $\lambda\in P^+$ and $\bos\in\bz_+^I$ satisfy
 \begin{enumerate}[{\rm(i)}]
\item any non--zero element of $\Hom_{\lie g\otimes A}(W_A(\lambda), W_A(\bos))$ is injective,
 \item   the map $\tau_\bos: W_A(\mu_\bos)\to W_A(\bos)$ is injective.
 \end{enumerate}
     Then any non--zero element of $\Hom_{\lie g\otimes A}(W_A(\lambda), W_A(\mu_\bos))$ is injective. Moreover, if $s_i=0$ for all $i\notin I_0$, then
      $$\Hom_{\lie g\otimes A}(W_A(\lambda), W_A(\mu_\bos))=0, \qquad \lambda\ne \mu_\bos.$$ \end{lem}
     \begin{pf} Let $\eta\in\Hom_{\lie g\otimes A}(W_A(\lambda), W_A(\mu_\bos))$. If $\eta\ne 0$, then $\tau_\bos\cdot\eta\in\Hom_{\lie g\otimes A}(W_A(\lambda),
     W_A(\bos))$ is non--zero since $\tau_\bos$ is injective. Hence $\tau_\bos\cdot\eta$ is injective which forces $\eta$ to be injective. If we now assume that
     $\lambda\ne\mu_\bos$ and that $s_i=0$ if $i\notin I_0$, then it follows from Corollary \ref{class} that $\Hom_{\lie g\otimes A}(W_A(\lambda), W_A(\bos))=0$ and
     hence it follows that $\eta=0$ in this case.
     \end{pf}
\begin{rem}\label{red} Using Theorem \ref{fundg},  we see that if $\lie g$ is of type $B_n$ or $D_n$ with $n\ge 6$ and $i=4$, then $(\bw^{\omega_4}_AA/\mathfrak I)^{\lie n^+\otimes A}_{\omega_2}\ne 0$ or equivalently $$\Hom_{\lie g\otimes A}(\bw^{\omega_2}_AA/\mathfrak I,\bw^{\omega_4}_AA/\mathfrak I)\ne 0.$$ Using Proposition \ref{glfundreal} and \eqref{inc} we get $$\Hom_{\lie g\otimes A}(W_A(\omega_2),W_A(\omega_4))\ne 0,$$ which in particular proves that the last assertion of the above Lemma and hence
Theorem \ref{thm2} fail in this case.
\end{rem}

       \subsection{}
       From now on, we shall assume that $\bos\in\bz_+^I$ is such that $s_i=0$ if $i\notin I_0$. By Corollary~\ref{class}, we see that $\Hom_{\lie g\otimes A}
       (W_A(\lambda), W_A(\mu_\bos))=0 $ if $\lambda\ne\mu_\bos$  and the first condition of   Lemma~\ref{conseq1} is trivially satisfied. Hence Theorem \ref{thm2} will
       follow if we show that $\mu_\bos$ satisfies both conditions in Lemma \ref{conseq1}. By the discussion at the start of Section 5, we see that proving that
       $\mu_\bos$ satisfies the first condition is equivalent to proving that $\tau_\bos\boa$ is injective for all $\boa\in\ba_\bos$.
       In other words, Theorem~\ref{thm2} follows if
       we establish the following.
        \begin{prop}\label{injectivea} Let $\bos\in\bz_+^I$ be such that $s_i=0$ if $i\notin I_0$.  For all $\boa\in\ba_\bos$, the canonical map $\tau_\mu\boa:
       W_A(\mu_\bos)\to W_A(\bos)$ given by extending
         $w_{\mu_\bos}\to w_\bos\boa$ is injective, in the following cases:
         \begin{enumerate}[{\rm(i)}]\item $A=\mathcal{R}_{0,1}$ or $\mathcal{R}_{1,0}$,
         \item $A=\cal R_{k,\ell}$, $\lie g=\lie{sl}_{n+1}$ and $\bos=(s,0,
         \dots,0)\in\bz_+^I$, $s>0$.\end{enumerate}
        \end{prop}
     The rest of the section is devoted to proving the proposition.

\subsection{} We begin by proving the following Lemma. \begin{lem} Let $A$ be a finitely generated integral domain. Let $\bos\in\bz_+^I$ be such that $s_i=0$ if $i\notin I_0$. Then $\tau_{\mu_\bos}\boa: W_A(\mu_\bos)\to W_A(\bos)$ is injective for $\boa\in\ba_\bos\setminus\left\{0\right\}$ if and only if $\tau_{\mu_\bos}$ is injective.\end{lem}
\begin{pf} Consider the map $\rho_\boa:W_A(\bos)\to W_A(\bos)$ given by $$\rho_\boa(w)=w\boa,\qquad w\in W_A(\bos).$$ This is clearly a map of $(\lie g\otimes A,\ba_\bos)$--bimodules. Since $$W_A(\bos)_{\mu_\bos}=w_\bos\otimes\ba_\bos,$$ and $\ba_\bos$ is an integral domain, it follows that that the restriction of $\rho_\boa$ to $W_A(\bos)_{\mu_\bos}$ is injective, and so $$\ker\rho_\boa\cap W_A(\bos)_{\mu_\bos}=\{0\}.$$ Since $\wt W_A(\bos)\subset\mu_\bos-Q^+$, it follows that if $\ker\rho_\boa$ is non--zero, there must exist $w'\in\ker\rho_\boa$ with $$(\lie n^+\otimes A)w'=0.$$ But this is impossible by Corollary~\ref{class} and hence $\ker\rho_\boa=0$. Since $\tau_{\mu_\bos}\boa=\rho_\boa\tau_{\mu_\bos}$, the Lemma follows.
\end{pf}

  \subsection{} We now prove that $\tau_{\mu_\bos}$ is injective. This was proved in \cite{CPweyl} for $\lie g=\lie{sl}_2$ and $A=\mathcal{R}_{1,0}$  and in \cite{FL} for $\lie g=\lie{sl}_{n+1}$, $\bos=(s,0,\dots,0)\in\bz_+^n$, $s>0$ and for any finitely generated integral domain $A$.

 Since $$\tau_\bos W_A(\mu_\bos)_{\mu_{\bos}}\cong_{\ba_{\mu_\bos}}(\bu(\lie g\otimes A)w_\bos)_{\mu_\bos}\cong_{\ba_{\mu_\bos}}\ba_{\mu_\bos},$$ the following proposition completes the proof of Proposition \ref{injectivea}.

  \begin{prop}\label{injectiveb} Let   $\mu\in P^+$ and let $\bpi:W_A(\mu)\to W$ be a surjective map of $(\lie g\otimes A,\ba_\mu)$-bi-modules such that the restriction of $\bpi$ to $ W_A(\mu)_\mu$ is an isomorphism of right $\ba_\mu$--modules. If $A=\mathcal{R}_{0,1}$ or $\mathcal{R}_{1,0}$ and $\mu=\sum_{i\in I_0}s_i\omega_i$, then $\bpi$ is an isomorphism.
  \end{prop}

\subsection{} \label{subs-localweyl-struct} Assume from now on that $A$ is either $\mathcal{R}_{0,1}$ or $\mathcal{R}_{1,0}$.
The following is well--known.

\begin{prop} For all $r\in\bz_+$, the ring $(\mathcal{R}_{0,1}^{\otimes r})^{S_r}$ is isomorphic to $\mathcal R_{0,r}$ and $(\mathcal{R}_{1,0}^{\otimes r})^{S_r}$ is isomorphic to $\bc[t_1,t_2,\dots, t_r, t_r^{-1}]$.
\end{prop}

The proposition implies that $\Max \ba_\lambda$ is an irreducible variety.
Given~$\lambda\in P^+$, define
$\mathcal D_\lambda\subset\Max\ba_\lambda$ by: $\bi\in\cal D_\lambda$ if and only if the $S_{r_\lambda}$--orbit of $\sym_\lambda\bi$ is of maximal size,
i.e, $\sym_\lambda\bi$ is the $S_{r_\lambda}$--orbit of
$((t-a_{1,1}),\dots,(t-a_{1,r_1}),\dots,(t-a_{n,1}),\dots,(t-a_{n,r_n}))\in(\Max A)^{\times r_\lambda}$ for some $a_{i,r}\in\bc$ (respectively $a_{i,r}\in\bc^\times$) with
$a_{i,r}\not=a_{j,s}$ if~$(i,r)\not=(j,s)$. { The set of such orbits is clearly Zariski open in $\Max \mathbb A_\lambda$. Since
$\sym_\lambda$ induces an isomorphism of algebraic varieties $\Max\mathbb A_\lambda\to\Max\ba_\lambda$, we conclude that
$\cal D_\lambda$ is Zariski open, hence is dense in $\Max\ba_\lambda$.}
Therefore, given any non--zero $\boa\in\ba_\lambda$ there exists $\bi\in\cal D_\lambda$ with $\boa\notin\bi$.

\subsection{}  We shall need the following theorem.
\begin{thm}\label{localweyl-struct} Let $A=\mathcal{R}_{0,1}$ or $\mathcal{R}_{1,0}$ and let $\lambda=\sum_{i\in I}r_i\omega_i\in P^+$.
 \begin{enumerate}[{\rm(i)}]
 \item\label{localweyl-struct.i} The right $\ba_\lambda$--module $W_A(\lambda)$ is free of rank~$d_\lambda$, where $$d_\lambda=\prod_{i\in I}(\dim \bw^{\omega_i}_A(A/\mathfrak I))^{r_i},$$ for
any $\mathfrak I\in\Max A$.
     \item\label{localweyl-struct.ii}
     Let $\bi\in \mathcal D_\lambda$.
     Then
     $$\bw^\lambda_A(\ba_\lambda/\bi)\cong\bigotimes_{i\in I}\bigotimes_{r=1}^{r_i}\bw^{\omega_i}_A (A/\mathfrak I_{i,r}),$$ where $\mathfrak I_{i,r}\in\Max A$ is the ideal
     generated by $(t-a_{i,r})$. If, in addition, we have $r_i=0$ for $i\notin I_0$, then $\bw^\lambda_A(\ba_\lambda/\bi)$ is an irreducible $\lie g\otimes A$--
     module.
\end{enumerate}
\end{thm}
Part~\eqref{localweyl-struct.i} of the Theorem was proved in \cite{CPweyl} for $\lie{sl}_2$, in \cite{CL} for $\lie{sl}_{r+1}$ and in \cite{FoL} for algebras of type $A,D,E$.
The general case can be deduced from the quantum case, using results of \cite{BN, Ka,N}. Part~\eqref{localweyl-struct.ii} of the Theorem was proved in \cite{CPweyl} in a different language and in \cite{CFK} in the language of this paper.

\subsection{}
\begin{pf}[Proof of Proposition~\ref{injectiveb}]

Let $\{w_s\}_{1\le s\le d_\mu}$ be an $\ba_\mu$-basis of~$W_A(\mu)$ (cf. Theorem~\ref{localweyl-struct}\eqref{localweyl-struct.i}). Then for all $\bi\in\operatorname{max}\ba_\mu$, $\{w_s\tensor 1\}_{1\le s\le d_\mu}$ is a
$\bc$-basis of $\bw^\mu_A\ba_\mu/\bi$.
Suppose that $w\in\ker\bpi$ and write $$w=\sum_{s=1}^{d_\mu}w_s\boa_s,\qquad\boa_s\in\ba_\mu.$$
If $w\ne 0$, let $\boa$ be the product of the non--zero elements of the set $\{\boa_s:1\le s\le d_\mu\}$.
Since $\ba_\mu$ is an integral domain we see that $\boa\ne 0$. By the discussion in Section~\ref{subs-localweyl-struct}  we can choose $\bi\in\cal D_\mu $ with $\boa\notin\bi$. Then $\boa_s\not=0$ implies that $\boa_s\notin \bi$ and hence
$\bar w:=w\tensor 1=\sum_{s=1}^{d_\mu} w_s\tensor \bar\boa_s\not=0$, where $\bar\boa_s$ is the canonical image of~$\boa_s$ in~$\ba_\mu/\bi$.
Notice that Theorem~\ref{localweyl-struct}\eqref{localweyl-struct.ii} implies
 that $\bw_A^\mu (\ba_\mu/\bi)$ is a simple $\lie g\tensor A$-module.

Since $\bpi$ is surjective, $W$ is generated by $\bpi(w_\mu)$.
 Setting $W'=\bpi(W_A(\mu)\bi)$, we see that $$W'_\mu=\bpi((W_A(\mu)\bi)_\mu)=\bpi(w_\mu)\bi.$$
 In particular, this proves that $\bpi(w_\mu)\notin W'$, hence $W'$ is a proper submodule of $W$ and $$(W/W')\bi=0.$$ This implies that $\bpi$ induces a well--defined non--zero surjective homomorphism of $\lie g\otimes A$--modules $\bar\bpi: \bw^\mu_A(\ba_\mu/\bi)\to W/W'\to 0$. In fact since $\bw_A^\mu (\ba_\mu/\bi)$ is simple, we see that $\bar\bpi$ is an isomorphism. But now we have $$0=\bpi(w)=\bar\bpi(\bar w),$$ forcing $\bar w=0$ which is a contradiction caused by our assumption that $w\ne 0$.

The proof of Proposition \ref{injectiveb} is complete.
\end{pf}

\end{document}